\documentclass[12pt]{article}
\usepackage{amscd,amsmath,amssymb,amsfonts, dsfont}
\usepackage[nohug,PostScript=dvips]{diagrams}
\usepackage[cmtip, all]{xy}
\usepackage{pstricks}
\usepackage{color}
\usepackage{remark,theorem}
\font\lcirc=lcircle10 %scaled 800; from Paul Vojta
\def\cxdot{{\hskip2.5pt\raise3.8pt\hbox{\lcirc\char113}}}
\def\cx{{\textstyle{\cdot}}}

\long\def\ignore#1{}

\newcommand{\openbox}{\leavevmode
  \hbox to.77778em{%
  \hfil\vrule
  \vbox to.675em{\hrule width.6em\vfil\hrule}%
  \vrule\hfil}}
\newcommand{\qedsymbol}{\openbox}
\DeclareRobustCommand{\qed}{%
  \ifmmode % if math mode, assume display: omit penalty etc.
  \else \leavevmode\unskip\penalty9999 \hbox{}\nobreak\hfill
  \fi
  \quad\hbox{\qedsymbol}}

\newcommand{\ie}{\textit{i.e.}}

\newcommand{\proofname}{{ Proof:}}
\newenvironment{proof}[1][\proofname]{\par
 \normalfont
\trivlist
  \item[\hskip\labelsep\itshape
    #1
]\ignorespaces
}{%
 \qed\endtrivlist}
\def\oh#1{{\cal O}_{#1}}

\def\uar#1{\mathop{#1}\limits_{\raise3pt\hbox{$\to$}}}

\def\boxit#1{\vbox{\hrule\hbox{\vrule\kern1pt
    \vbox{\kern1pt#1\kern1pt}\kern1pt\vrule}\hrule}}

\def\varprojlim{\mathop{\vtop{\ialign{$##$\cr
 \hfil{\fam0lim}\hfil\cr\noalign{\nointerlineskip}%
 {\leftarrow}\mkern-6mu\cleaders\hbox{$\mkern-2mu{-}\mkern-2mu$}\hfill
 \mkern-6mu{-}\cr\noalign{\nointerlineskip\kern-.2326ex}\cr}}}}
\def\varinjlim{\mathop{\vtop{\ialign{$##$\cr
 \hfil{\fam0lim}\hfil\cr\noalign{\nointerlineskip}%
 {-}\mkern-6mu\cleaders\hbox{$\mkern-2mu{-}\mkern-2mu$}\hfill
 \mkern-6mu{\to}\cr\noalign{\nointerlineskip\kern-.2326ex}\cr}}}\nolimits}

\def\uar#1{\vec {#1}}

 \newcommand{\bc}{{\bf C}}
 \newcommand{\bq}{{\bf Q}}
 
 \newcommand{\bz}{{\bf Z}}

\newcommand{\et}{{\acute et}}

\newcommand{\spec}{\mathop{\rm Spec}\nolimits}

\newcommand{\im}{\mathop{\rm Im}\nolimits}

\newcommand{\Hom}{\mathop{\rm Hom}\nolimits}

\newcommand{\Gal}{\mathop{\rm Gal}\nolimits}

\newcommand{\Pic}{{\rm Pic}}

\newcommand{\ot}{\otimes}

\def\odiagram#1{
  \def\normalbaselines{\baselineskip20pt\lineskip3pt \lineskiplimit3pt }
   \matrix{#1}}
 \def\ldiagram#1{$$\odiagram#1 % display math mode hack
      \refstepcounter{eqncounter}\
      \edef\@currentlabel{\p@equation
 \thetheorem.\theeqncounter}} 
 \def\endldiagram{\eqno  \@eqnnum$$\global\@ignoretrue}
\def\ldiagram#1{$$\odiagram#1 % display math mode hack
      \refstepcounter{eqncounter}\
      \edef\@currentlabel{\p@equation
 \thetheorem.\theeqncounter}} 
 \def\endldiagram{\eqno  \@eqnnum$$\global\@ignoretrue}

\newtheorem{theorem}{Theorem}

\newtheorem{question}{Question}

\newtheorem{corollary}[theorem]{Corollary}
\newtheorem{proposition}[theorem]{Proposition}
\newtheorem{lemma}[theorem]{Lemma}
\newtheorem{claim}[theorem]{Claim}

\newremark{example}[theorem]{Example}
\newremark{remark}[theorem]{Remark}
\theorembodyfont{\normalfont\slshape}
\newcommand{\cH}{\mathcal{H}}
\newcommand{\cA}{\mathcal{A}}

\newcommand{\Gamm}{\Gamma}
\renewcommand{\bq}{\mathbb{Q}}
\title{Hodge cohomology of invertible sheaves}
\author{H\'el\`ene Esnault and Arthur Ogus}
\begin{document}
\maketitle
\section{Introduction}

Let $k$ be an algebraically closed field and let $X/k$
be a smooth projective  connected $k$-scheme.
Let $L$ be an invertible sheaf on $X$, and
for each integer $m$, let
$$H^m_{Hdg}(X/k,L) := \bigoplus_{a+b = m} H^b(X, L\otimes \Omega^a_{X/k}).$$
We wish to study how the dimensions of  the $k$-vector spaces
 $H^m_{Hdg}(X/k,L)$ and  $ { H^b(X, L\otimes \Omega^a_{X/k})}$ vary
with $L$.  
For example, if $k$ has characteristic zero,
Green and Lazarsefeld \cite{GL} proved that for 
 given $i,j,m$,  the subloci
$$\{ L\in {\rm Pic}^0(X): {\rm dim}\ H^i(X, \Omega_X^j\otimes L)\ge m\} $$
of $\Pic^0(X)$  are translates of abelian subvarieties, and   Simpson \cite{S} showed
that they in fact  are translates by torsions points.  
Both these papers use analytic methods,
but  Pink and Roessler \cite{PR} obtained the same results purely algebraically,
using  the technique of 
 mod $p$ reduction  and the decomposition theorem of Deligne-Illusie. 
A key point of their proof is the fact that
if that if $L^n \cong \oh X$ for some positive integer $n$, 
 then for all natural numbers $a$ with $(a,n)=1$ one has
\begin{gather}\label{eq2} \dim H^m_{Hdg}(X/k,L)=\dim H^m_{Hdg}(X/k,L^a) \end{gather}
 (\cite[Proposition~3.5]{PR}).
They conjecture that equation~\ref{eq2} remains true 
in  characteristic $p>0$ if $X/k$ lifts to $W_2(k)$ and has dimension $\le p$.
The purpose of this note is to discuss a few aspects of this  conjecture and some variants.

Our main result (see Theorem~\ref{np.t}) says 
that the conjecture is true if $n = p$ and 
 $X$ is ordinary in the sense of Bloch-Kato~ 
\cite[Definition~7.2]{BK}.
 We also explain in section \ref{mv.s} some motivic
variants of  (1) and, in particular  in Proposition~\ref{prop:chow}, 
a proof (due to Pink and Roessler) of the characteristic zero case of 
 (\ref{eq2}), using the language of Grothendieck Chow motives. 
See \cite[9.3]{KA} for a discussion of a related problem using
similar techniques.
We should remark that there  are also some log versions of these questions,
which we will not make explicit.
\\[.2cm]
{\it Acknowledgements:} We thank D. Roessler for explaining to us his and R. Pink's  analytic proof of equation~\ref{eq2}. We thank the referee for very useful, accurate and  friendly remarks which helped us improving the exposition of this note.

\section{Motivic  variants}\label{mv.s}

\begin{question}\label{qt.q}
Let $X$ be a smooth projective connected variety defined over an
 algebraically closed field $k$. Let $L$ be an invertible sheaf on $X$ and $n$
a positive integer such that
 $L^n \cong \oh X$.  Is 
$$\dim H^m_{Hdg}(X/k,L^i) = \dim H^m_{Hdg}(X/k,L)$$
for every $i$ relatively prime to $n$?
\end{question}

Let us explain how this question can be given a motivic
interpretation.
 We refer to \cite{Sch} for the definition of Grothendieck's Chow motives over  a  field $k$. In particular, objects are triples $(Y,p,n)$ 
where $Y$ is a smooth projective variety over $k$, 
 $p$ is an element $CH^{{\rm dim}(Y)} (Y\times_k Y)\otimes \mathbb{Q}$,
 the rational Chow group of ${\rm dim}(Y)$-cycles, which, as a correspondence, 
is an idempotent, and $n$ is a natural number. 

Let $\pi: Y\to X$ be a principal bundle  under a $k$-group scheme  $\mu$,
where $X$ and $Y$ are smooth and projective over $k$.  
Recall that this means that there is a $k$-group scheme action $\mu\times_k Y\to Y$ with the property that one has an isomorphism
$$ (\xi,y) \mapsto (y,\xi y) :
\mu \times_k Y \cong Y \times_X Y \subseteq Y\times_k Y.$$
Thus a point $\xi \in \mu(k)$ defines a  closed subset
$\Gamma_\xi$ of $Y \times_k Y$,  the graph of the 
endomorphism of $Y$ defined by $\xi$. The map 
$\xi \mapsto \Gamma_\xi$  extends  uniquely to a map of ${\mathbb Q}$-vector spaces
$$ \Gamma: {\mathbb Q}[\mu(k)] \to  CH^{{\rm dim}(Y)}(Y\times_k Y) \otimes {\mathbb Q}.$$
Here ${\mathbb Q}[\mu(k)] $ is the ${\mathbb Q}$-group algebra,
so the product structure is induced by the product of $k$-roots of unity. We can think of 
$CH^{{\rm dim}(Y)}(Y\times_k Y) \otimes {\mathbb Q}$ as a ${\mathbb Q}$-algebra of
 correspondences acting on $CH^*(Y)\otimes {\mathbb Q}$, where for $\beta \in {CH^s(Y)\otimes {\mathbb Q}}, \gamma\in  CH^{{\rm dim}(Y)}(Y\times_k Y) \otimes {\mathbb Q}$, one defines as usual 
$$\gamma\cdot \beta :=(p_2)_*(\gamma\cup p_1^*\beta).$$
Then the  map $\Gamma$ is easily seen to be   compatible with composition, as on closed points $y \in Y$ one has $\Gamma_\xi(y)=\xi\cdot y$. 
In particular if $\xi \in \bq[\mu]$ is idempotent
in the group ring $\bq[\mu(k)]$, then $\Gamma_\xi\cong Y\times \xi$ is idempotent as a correspondence. 
  In this case we let  $Y_\xi$  be the Grothendieck Chow motive
$(Y, \xi, 0)$.

Let $L$ be an $n$-torsion invertible sheaf on  smooth irreducible projective scheme $X/k$. Recall that the choice of 
an $\oh X$-isomorphism $L^n\xrightarrow{\alpha} \mathcal{O}_X$ 
defines an $\oh X$-algebra structure on  
\begin{equation}\label{adef.e}
\mathcal{A} := \bigoplus_{i= 0}^{n-1} L^i
\end{equation}
via the tensor product $L^i\times L^j\to L^{i}\otimes_{\mathcal{O}_X} L^j=L^{i+j}$ for $i+j<n$
and its composition with the isomorphism 
$L^i\times L^j\to L^{i}\otimes_{\mathcal{O}_X} L^j=L^{i+j} \xrightarrow{\alpha^{-1}} L^{i+j-n}$
for  $0\le i+j-n$. 
 Then the corresponding $X$-scheme
$\pi: Y:= \spec_X \cA \to X$
is a torsor under the group scheme  $\mu_n$ of $n$th roots of unity.
 Indeed, locally Zariski on $X$, $\cA\cong \mathcal{O}_X[t]/(t^n-u) $
for a local unit $u$, the $\mu_n$-action is defined by $\cA\to 
\cA\otimes \mathbb{Q}[\zeta]/(\zeta^n-1), \ t\mapsto t\zeta$,
 and the torsor structure is given by $
 \cA \otimes   \mathbb{Q}[\zeta]/(\zeta^n-1)\cong
\cA\otimes_{  \mathcal{O}_X} \cA, \ (t,\zeta)\mapsto (t, t\zeta) $.  
This construction defines an equivalence between the category
of pairs $(L,\alpha)$ and the category of $\mu_n$-torsors over $X$. 
  Assuming now that $n$ is  invertible in $k$, $\mu_n$ is \'etale,
hence $\pi$ is \'etale and $Y$ is smooth  and projective over $k$.        
Note  the character group $X_n := \Hom(\mu_n,{\bf G}_m)$
is cyclic of order $n$ with a canonical generator
(namely, the inclusion $\mu_n \to {\bf G}_m$).
By construction, the direct sum decomposition~(\ref{adef.e}) of $\cA$
corresponds exactly to its eigenspace decomposition
according to the characters of $\mu_n$.  

 We can now apply the general construction of motives to this situation.
Since $\mu_n$ is \'etale over the algebraically closed field  $k$,
it is completely determined by the finite group  $ \Gamma :=\mu_n(k)$,
which is cyclic of order $n$.  The group algebra
$\bq[\Gamm]$ is a finite separable
algebra over $\bq$, hence is a product of fields:
$$\bq[\Gamm] = \prod E_e.$$  
Here $E_e=\bq[T]/(\Phi_e(T))=\bq(\xi_e)$, where $e$ is a divisor of $n$,  $\Phi_e(T)$ is
the cyclotomic polynomial, and $\xi_e$ is a primitive $e$th root of unity. 
   There
is an (indecomposable) idempotent $e$ corresponding
to each of these fields, and
for each $e$  we find a Chow motive $Y_e$.

The indecomposable idempotents of $\bq[\Gamm]$ 
can also be thought of as points of the spectrum $T$
of $\bq[\Gamma]$.  
If $K$ is  a sufficiently large extension of $\bq$, then
\begin{equation}
  \label{eq:5}
T(K) = \Hom_{\rm Alg}(\bq[\Gamma ], K) = \Hom_{\rm Gr}(\Gamm,K^*),
\end{equation}
\begin{equation}
  \label{eq:6}
\mbox{and     }K\ot \bq[\Gamm] \cong K[\Gamm] \cong K^{T(K)}.
  \end{equation}

Thus $T(K)$ can be identified with the character group $X_n$ of $\Gamma$,
and is canonically isomorphic to $\bz/n\bz$, 
with canonical generator the inclusion $\Gamma \subseteq k$. 
Suppose that $K/\bq$ is Galois.  Then
${\rm Gal}(K/\bq)$ acts on $T(K)$, and  the points of $T$
correspond to the $\Gal(K/\bq)$-orbits.  
By the theory of cyclotomic extensions of $\bq$, this action
factors through a surjective map
$${\rm Gal}(K/\bq) \to (\bz/n\bz)^*$$
and the usual action of $(\bz/n\bz)^*$ on $\bz/n\bz$
by multiplication.
Thus the orbits correspond precisely to 
the divisors $d$ of $n$; we shall associate
to each orbit $S$  the index $d$   of the subgroup of $\bz/n\bz$
generated by any element of $S$.  (Note that in fact
the image of $d$ in $\bz/n\bz$ belongs to $S$.)
 We shall  thus identify
the indecomposable idempotents of $\bq[\Gamm]$ and the divisors of
$n$. 

Let us suppose that $k = \bc$.  Then we can consider
the Betti cohomologies of $X$ and $Y$, and in particular
the group algebra $\bq[\Gamm]$ operates on
 $H^m(Y,\bq)$.  We can thus view
$H^m(Y,\bq)$ as a $\bq[\Gamm]$-module,
which corresponds to  a coherent sheaf $\tilde H^m(Y,\bq)$
on  $T$.  If $e$ is an idempotent of $\bq[\Gamm]$, then 
$H^m(Y_e,\bq)$ is the image of  the action  of $e$ 
on $H^m(Y,\bq)$, or equivalently, it is the stalk
of the sheaf $\tilde H^m(Y,\bq)$ at the point of $T$ corresponding
to $e$, or equivalently, it is
$H^m(Y,\bq)\ot E_e$ where the tensor product is taken
over $\bq[\Gamm]$.  If $K$ is a sufficiently large 
field as above, then equation (\ref{eq:6}) induces an isomorphism
of $K$-vectors spaces:
$$H^m(Y_e,\bq)  \ot_{\mathbb{Q}}K \cong  
\bigoplus \{ H^m(Y,K)_t: {t \in T^e(K)} \},$$
where here $T^e(K)$ means the set of points of $T(K)$
in the Galois orbit corresponding to $e$, and 
$H^m(Y,K)_t$ means the $t$-eigenspace of the action of $\Gamm$
on $ H^m(Y,\bq){\ot_{\mathbb{Q}} } K $.    The de Rham and Hodge cohomologies
of $Y_e$ are defined in the same way: they are the 
images of the actions of the idempotent $e$
acting on the $k$-vector spaces $H_{DR}(Y/k)$
and $H_{Hdg}(Y/k)$.

The following result is due to Pink and Roessler.  Their
article \cite{PR} contains a proof using reduction
modulo $p$ techniques and the results of \cite{DI}; the following
analytic argument is based on oral communications with them.

\begin{proposition}\label{prop:chow}
The answer to question 1 is affirmative if $k$ is a field of
characteristic zero.
\end{proposition}
\begin{proof}
As both sides of the equality in Question 1 satisfy base
change with respect to field extensions, 
we may assume that  $k = \bc$. 
Let $i \to t_i$ denote the isomorphism
 $\bz/n\bz \cong T(\bc)$. 
For each divisor $e$ of $n$ there is a corresponding
idempotent $e$  of $\bq[\Gamma] \subseteq K[\Gamma]$,
 the sum over all $i$ such that  $t_i  \in 
T^e(\bc)$.  Consider the Hodge
cohomology of the motive $Y_e$:
\begin{eqnarray*}
H^m_{Hdg}(Y_e/\bc) := H^m_{Hdg}(Y/\bc)\ot_{\bq[\Gamma]}{  E_e} 
    &\cong& H^m_{Hdg}(Y/\bc)\ot_{\bc[\Gamma]}{  (\bc\ot E_e)}. \\
    & \cong & { \bigoplus \{ H^m_{Hdg}(Y/\bc)_i : {i \in T^e(k)} \}}.
\end{eqnarray*}
Since $ \pi \colon Y \to X$ is  finite and  \'etale, 
\begin{eqnarray*}
H^b(Y,\Omega^a_{Y/\bc}) \cong H^b(X,\pi_*{ \pi^* \Omega^a_{X/\bc}}) &\cong&
H^b(X, \Omega^a_{X/\bc} \ot \pi_*\oh Y)\\
&\cong&  \bigoplus \{ H^b(X,\Omega^a_{X/\bc} \ot L^i) : i \in \bz/n\bz\}.
\end{eqnarray*}
Thus  
$$H_{Hdg}^m(Y/\bc) \cong  \bigoplus \{  H^m_{Hdg}(X, L^i): i \in \bz/n\bz\}, $$
and hence  from the explicit
description of the action of $\mu_n$ on $\cA$ above
it follows that 
$$H^m_{Hdg}(Y_e/\bc) =   \bigoplus\{ H^m_{Hdg}(X, L^i): i \in T_e(\bc {)} \}.$$ 
The Hodge decomposition theorem for $Y$ provides us with an isomorphism:
$$H^m_{Hdg}(Y/\bc) \cong \bc\otimes H^m(Y,\bq),$$
compatible with the action of $\bq[\Gamma]$.  This gives us,
for each idempotent $e$, an isomorphism of $\bc\ot E_e$-modules.
$$H^m_{Hdg}(Y_e/\bc) \cong \bc\otimes H^m(Y_e,\bq).$$
  The action on $\bc\otimes H^m(Y_e,\bq)$
on the right just comes from the action of $E_e$ on $H^m(Y_e,\bq)$
by extension of scalars.  Since $E_e$ is a field, $H^m(Y_e,\bq)$
is free as an $E_e$-module, and hence
the $\bc\ot E_e$-module $H^m_{Hdg}(Y_e/\bc)$ is also free.
It follows that its rank is the same at all the points $t \in T^e(\bc)$,
affirming Question~\ref{qt.q}.
\end{proof}

Let us now formulate an  analog of Question 1 for the $\ell$-adic
and crystalline realizations of the motive $Y_e$ in characteristic $p$. 

\begin{question}\label{q2}
  Suppose that $k$ is an algebraically closed field of characteristic $p$
and $(n,p) = 1$.  Let $\ell$ be a prime different from $p$,
let $e$ be a divisor of $n$, and let $E_e$ be the corresponding
factor of $\bq[\Gamma]$.
Is it true that each $H^m(Y_e,\bq_\ell)$ is a free $\bq_\ell\ot E_e$-module?
And is it true that $H^m_{cris}(Y_e/W)\ot \bq$ is a free
$W \ot E_e$-module, where $W:= W(k)$? 
\end{question}
If $K$ is an extension of $\bq_\ell$ (resp. of $W(k)$)  which
contains a primitive $n$th root of unity,
then as above we have a eigenspace decompositions:
\begin{eqnarray}
  \label{eq:7}
  K\ot H^m(Y_\et,\bq_\ell)  &\cong&  \bigoplus \{ H^m(Y_\et,K)_t: t\in T(K)\}\\
  K\ot H^m(Y_{cris}/W(k))  &\cong&  \bigoplus \{ H^m(Y_{cris},K)_t: t\in T(K)\},
\end{eqnarray}
and this question asks whether the $K$-dimension of the $t$-eigenspace
is constant over the orbits $T_e(K) \subseteq T(K)$. 

We show in the sequel that the question has a positive answer.

Suppose first that $X/k$ lifts to characteristic zero, \ie,
that there exists a complete discrete valuation ring $V$
with residue field $k$ and fraction field of characteristic zero
and a smooth proper $\tilde X/V$ whose special fiber is $X/k$.
Let $X_m$ be the closed subscheme of $\tilde X$ defined by
$\pi^{m+1}$, where $\pi$ is a uniformizing parameter of $V$. 
Choose a trivialization $\alpha$ of $L^n$.
It follows from Theorem 18.1.2 of \cite{EGAIV4}
 that the \'etale $\mu_n$-torsor $Y$
on $X$  corresponding to $(L,\alpha)$ lifts to $X_m$, uniquely up
to a unique isomorphism, and hence that the same is true for
$(L,\alpha)$.  This fact can also be seen by chasing the  exact
sequences of cohomology corresponding to the
commutative diagram of exact sequences in the \'etale topology
\begin{gather}
\xymatrix{ & & 0\ar[d]  & 0\ar[d] \\
 & & \ar[d]_{a\mapsto 1+\pi^ma}\mathcal{O}_X \ar[r]^{n \ \cong}& \mathcal{O}_X \ar[d]^{a\mapsto 1+\pi^ma}\\ 
1\ar[r] &  \ar[d]_{=} \mu_n \ar[r] & \ar[d]  \mathcal{O}^\times_{X_{m}}  \ar[r]^{n} & \ar[d] \mathcal{O}^\times_{X_{m}} \ar[r] & 1\\
1\ar[r] &  \mu_n \ar[r] & \ar[d] \mathcal{O}^\times_{X_{m-1}}  \ar[r]^{n} & \ar[d]  \mathcal{O}^\times_{X_{m-1}} \ar[r] & 1 \\
& & 1 & 1} 
\end{gather}
 By Grothendieck's fundamental theorem for proper morphisms, 
it follows  that  $(L, \alpha)$ and $Y$  lift to 
 $(\tilde{L}, \tilde{\alpha})$ and $\tilde Y$ on $\tilde X$. 
Then by the \'etale to Betti and Betti to crystalline
comparison theorems, we see
that under the lifting assumption, 
the answer to Question \ref{q2} is affirmative.

In fact, the lifting hypothesis is superfluous,
but this takes a bit more work.

\begin{claim}
  The answer to Question \ref{q2} is affirmative.
\end{claim}
\begin{proof}
It is trivially true that $H^m(Y_e,\bq_\ell)$ is free over $\bq_\ell \ot E_e$
if  $\bq_\ell \ot E_e$ is a field. If $(\ell,n) = 1$, 
 this is the case if and only
if $(\bz/e\bz)^*$  is cyclic and generated by $\ell$. 
More generally, assuming $\ell$ is relatively prime to $n$,
 there is a decomposition of $\bq_\ell[\Gamma]$ into a product of fields
$\bq_\ell[\Gamma] \cong  \prod E_{\ell,e}$, where now $e$
ranges over the orbits of $\bz/n\bz$ under the action
of the cyclic subgroup of $(\bz/n\bz)^*$ generated by $\ell$.  
This is indeed the unramified lift of the decomposition
 of $A=\mathbb{F}_\ell[\Gamma]$ into a product of finite extensions of $\mathbb{F}_\ell$, 
corresponding to the orbits of Frobenius 
on the geometric points of $A$.
This shows   at least that the dimension of
$H^m(Y,K)_t$ in (\ref{eq:7}) is, as a function of $t$, constant
over the $\ell$-orbits.  

For the general statement, let $K$ be an algebraically closed field containing
$\bq_\ell$ for all  primes $\ell \neq p$, and containing
$W(k)$.  For  $\ell \neq p$
let $V_\ell := H^m(Y_\et,\bq_\ell)\ot_{\bq\ell} K$, and 
let $V_p := H^m(Y_{cris},W(k))\ot_{W(k)} K$.
Then each $V_\ell$ is a finite-dimensional representation
of $\Gamma$, and the isomorphisms (\ref{eq:7})  and (6) are just its
decomposition as a direct sum of irreducible representations:
$$V_\ell \cong  \bigoplus \{  n_{\ell,i} V_i :i \in \bz/n\bz \},$$
 where $V_i = K$, with $\gamma \in \Gamma$ acting
by multiplication by $\gamma^i$.  
By \cite[Theorem~2.2)]{KM} (and \cite{Berth}, \cite{Gros} and \cite{Pe} for the existence of cycle classes 
in crystalline cohomology)
the trace of  any  $\gamma 
\in \Gamma$  
acting on $V_\ell$ is an integer independent of $\ell$, including $\ell = p$. 
 Since $\Gamma$ is a finite group,
it follows from the independence of characters that for 
each $i$,  $n_i :=n_{\ell,i}$
is  independent of $\ell$.  We saw above that
$n_{\ell,\ell i} = n_{\ell, i}$ if $(\ell, n) = 1$ and
$\ell \neq p$, so that in fact $n_{\ell i} = n_i$ 
for all $\ell \neq p$ with $(\ell, n) = 1$. 
Since the group $(\bz/n\bz)^*$ is generated by all such 
$\ell$, it follows that $n_i$ is indeed constant over the $\ell$-orbits.
\end{proof}

What does this tell us about Question 1?  If $(p,n) = 1$ and $k$
is algebraically closed, $W[\Gamma]$ is still semisimple,
and can be written canonically as a product of copies of $W$,
indexed by $i \in T(W) \cong \bz/n\bz$. 
For every $t \in T(W) \cong T(k)$,
we have an injective base change map from crystalline to de Rham cohomology:
$k\ot H^m(Y/W)_t \to H^m(Y/k)_t$.

\begin{question}\label{q3}
  In the above situation, is $H^q(Y/W)$ torsion free when $(p,n) = 1$?
\end{question}
If the answer is yes, then the maps $k\ot H^m(Y/W)_t \to H^m(Y/k)_t$
 are isomorphisms,
and this means that we can compute the dimensions
of the de Rham eigenspaces from the $\ell$-adic ones.
Assuming also that the Hodge to de Rham spectral sequence
of $Y/k$ degenerates at $E_1$, this should give an affirmative
answer to Question 1.  Note that if $X/k$ lifts mod $p^2$,
then  $Y/k$ lifts mod $p^2$ as well, and if the dimension is less than
or equal to $p$, the $E_1$-degeneration is true by \cite{DI}.

Of course, there is no reason for Question \ref{q3} to have
an affirmative answer in general.  Is there a reasonable
hypothesis on $X$ which guarantees it?  For example,
is it true if the crystalline cohomology of $X/W$
is torsion free?

\section{The $p$-torsion case in characteristic $p$}
Let us assume from now on that $k$ is a perfect field of  characteristic $p>0$.
In this case we can reduce question \ref{qt.q} to 
a question about connections, using the following construction
of \cite{DI}.  First let us recall some standard notations.
 Let $X'$ be the pull back of $X$ via the Frobenius of $k$,
let $\pi \colon X' \to X$ be the projection, and let
$F: X\to X'$ and $F_X \colon X \to X$ be the relative  and
 absolute Frobenius morphisms. Then $F_X^*L=L^p=F^*L'$, where $L':= \pi^*L$. 
 Then  $L^p=F^{-1}L'\otimes_{\mathcal{O}_{X'}}\mathcal{O}_{X}$ is endowed 
with the Frobenius descent connection $1\otimes d$, {\em i.e.} the unique connection spanned by its flat sections $L'$.  In general, for a given integrable connection $(E,\nabla)$, we set 
$$H^i_{DR}(X, (E,\nabla))={\mathbb H}^i(X/k, (\Omega^\cx_{X/k} \otimes E, \nabla)),$$
 and we use again the notation
$$H^i_{Hdg}(X/k, L)=\bigoplus_{a+b=i}H^b(X, \Omega^a_{X/k}\otimes L)$$
and write $h^i_{DR}$ and $h^i_{Hdg}$ for the respective  dimensions of these spaces.

\begin{proposition}\label{di.p}
  Let $L$ be an invertible sheaf on a smooth proper  scheme $X$   over   $k$ and let
$\nabla$ be the Frobenius descent connection on $L^p$.
Suppose that $X/k$ lifts to $W_2(k)$ and has dimension
at most $p$.
Then for every natural number $m$,
$$h^m_{DR}(X/k,(L^p,\nabla)) = h^m_{Hdg}(X/k,L).$$\end{proposition}
\begin{corollary} \label{cor}
Under the assumtpions of Proposition \ref{di.p},  if $L^p \cong \oh X$ and $\omega := \nabla(1)$,
then for any integer $a$,
$$h^m_{Hdg}(X/k,L^a) = h^m_{DR}(X/k, (\oh X,d+a\omega)).$$
\end{corollary}
\begin{remark}  If $p$ divides $a$, this just means the degeneration of the Hodge to de Rham spectral sequence for $(\oh X,d)$. \end{remark}

\begin{proof}
Let $Hdg^\cx_{X'/k}$ denote the Hodge complex of $X'/k$,
\ie, the direct sum $\oplus_i\Omega^i_{X'/k}[-i]$.
  Recall from \cite{DI} that { the lifting yields }
an isomorphism in the { bounded} derived category of $\oh {X'}$-modules:
$$Hdg^\cx_{X'/k} \cong F_*(\Omega^\cx_{X/k},{ d}).$$
Tensoring this isomorphism with $L' := \pi^*L$ and using the
projection formula for $F$, we find an isomorphism
$$Hdg^\cx_{X'/k}\ot L' \cong F_*(\Omega^\cx_{X/k}\ot L^p,\nabla).$$
  Hence 
$$  H^m_{Hdg}(X/k,L) \xrightarrow{F_k^* \ \cong} H^m_{Hdg}(X'/k,L') 
\xleftarrow{F_* \ \cong }  H^m_{DR}(X,(L^p,\nabla)).$$  This proves the proposition. If $L^p = \oh X$,
 the corresponding Frobenius descent connection
$\nabla$ on $\oh X$ is determined by $\omega_L := \nabla(1)$.
It follows from the tensor product rule for connections
that $\omega_{L^a} = a\omega_L$ for any integer $a$. 
\end{proof}

The corollary suggests the following question.

 \begin{question}\label{qomega.q}
   Let $\omega$ be a closed one-form on $X$ and let $c$ be a unit of $k$.
Is the dimension of $H_{DR}^m(X,(\oh X, d+c\omega))$ independent of $c$?
 \end{question}

\begin{remark}
Some properness is necessary, since the $p$-curvature
of $d_\omega:= d+ \omega$
can change from zero to non-zero as one multiplies by an
invertible constant.
If the $p$-curvature is non-zero, then the sheaf $\cH^0(\Omega_{X/k}^\cx,d_\omega)$
vanishes, and hence so does $H^0(X,(\Omega_{X/k}^\cx,d_\omega))$.  If the
$p$-curvature vanishes, then $\cH^0(\Omega_{X/k}^\cx,d_\omega)$
is an invertible sheaf $L$, which can have nontrivial sections
if $X$ is allowed to shrink. However, since by definition,  $L\subset \mathcal{O}_X$, it can have a global section
 on a proper  $X$ only if $L=\mathcal{O}_X$.  
\end{remark}

We can answer Question \ref{qomega.q}  under a strong hypothesis.

\begin{theorem}\label{np.t}
Suppose that $X/k$ is smooth, proper, and ordinary in the
sense of  Bloch and 
Kato \cite[Definition~7.2]{BK}:  $H^i(X,B^j_{X/k}) = 0$ for all $i,j$,
where 
$$B^j_{X/k} := Im\left (d \colon \Omega^{j-1}_{X/k} \to \Omega^j_{X/k}\right).$$
Then the answer to question \ref{qomega.q} is affirmative.
Hence if $X/k$ lifts to $W_2(k)$,  has dimension at most $p$,
and if $n = p$, the answer to Question $\ref{qt.q}$ is also affirmative.
\end{theorem}
We begin with the following lemmas.

 \begin{lemma}\label{bd.l}
Let $\omega$ be a closed one-form on $X$,
 and let 
$$d_\omega := d + \omega \wedge \quad \colon \quad\Omega^\cx_{X/k} \to \Omega^{\cx+1}_{X/k}. $$
Then the standard exterior derivative induces a morphism of complexes:
$$(\Omega_{X/k}^\cx,d_\omega) \rTo^\delta  (\Omega_{X/k}^\cx,d_\omega)[1]{\red .}$$ 
 \end{lemma}
 \begin{proof}
If $\alpha$ is a section of   $\Omega^q_{X/k}$,
\begin{eqnarray*}
  d d_\omega(\alpha) & = & d(d\alpha + \omega \wedge \alpha) \\
         & = & dd \alpha +d\omega \wedge \alpha - \omega \wedge d\alpha \\
  & = & -\omega \wedge d\alpha.
\end{eqnarray*}
Since the sign of the differential of the complex $(\Omega^\cx_{X/k},d_\omega)[1]$
is the negative of the sign of the differential of $(\Omega^\cx_{X/k},d_\omega)$,
\begin{eqnarray*}
  d_\omega d(\alpha) & = &-(d +\omega \wedge)(d\alpha) \\
  & = & -\omega \wedge d\alpha
\end{eqnarray*}
 \end{proof}

 \begin{lemma}\label{bz.l}
Let $Z^\cx := \ker(d)\subseteq (\Omega^\cx_{X/k},d_\omega)$ and 
$B^\cx := \im (d)[-1]\subseteq (\Omega^\cx_{X/k},d_\omega)$. 
Then for any $a \in k^*$, multiplication by $a^i$ in degree $i$ 
induces  isomorphisms
$$(Z^\cx, d_\omega) \rTo^{\lambda_a}  (Z^\cx,d_{a\omega})$$
$$(B^\cx, d_\omega) \rTo^{\lambda_a}  (B^\cx,d_{a\omega}){\red .}$$
 \end{lemma}
 \begin{proof}

It is clear that the boundary map $d_\omega$ on $Z^\cx$
and on $B^\cx$ is just wedge product with $\omega$.
\end{proof}

\begin{proof}[Proof of Theorem \ref{np.t}]
The morphism $\delta$ of Lemma \ref{bd.l} induces an exact sequence:
\begin{equation}
  \label{eq:8}
 0 \to (Z^\cx,d_\omega) \to  (\Omega^\cx_{X/k},d_\omega) \xrightarrow{\delta}
(B^\cx,d_\omega)[1] \to 0.
\end{equation}
As  $X/k$ is ordinary,    the $E_1$ term of the first
spectral sequence for $(B^\cx,d_\omega)$ is $E_1^{i,j} =
 H^j(X,B^i) = 0$,
and it follows that the hypercohomology of $(B^\cx,d_\omega)$ vanishes,
for every $\omega$.  Hence the natural map
$H^q(Z^\cx,d_\omega) \to  H^q(\Omega^\cx_{X/k},d_\omega)$ is an isomorphism.
Since 
the dimension of $H^q(Z^\cx,d_\omega)$ is unchanged when $\omega$ is multiplied
by a unit of $k$, the same is true of $ H^q(\Omega^\cx_{X/k},d_\omega)$.
This completes the proof of Theorem~\ref{np.t}.
\end{proof}

\begin{remark}
 A simple Riemann-Roch computation shows that on curves,   question 1  has a positive answer with no additional assumptions.
Indeed, if $L$ is a nontrivial torsion sheaf, then its degree is zero and it has no global sections.
It follows that $h^1(L) = g-1$.  Since the same is true for $L^{-1}$, $h^0(L\otimes \Omega^1_X) = h^1(L^{-1}) = g-1$,
and $h^1(L\otimes \Omega^1_X) = h^0(L^{-1}) = 0$. 
\end{remark}

\begin{remark}
  In the absence of the ordinarity hypothesis, one can
ask if the rank of the boundary map
$$\partial_\omega \colon H^{q+1}(B^\cx,\omega\wedge) \to H^{q+1}(Z^\cx,\omega\wedge) $$
of \eqref{eq:8} 
changes if $\omega$ is multiplied by a unit of $k$. 
To analyze this question, let 
$$c_\omega \colon (B^\cx,\omega\wedge) \to (Z^\cx,\omega\wedge)$$
be the morphism in the derived category $D(X', \mathcal{O}_{X'})$ defined
by the exact sequence \eqref{eq:8}, so that 
$\partial_\omega$ can be identified with $H^{q-1}(c_\omega)$.
Similarly, the exact sequence
$$  0 \to (Z^\cx,\omega\wedge) \to (\Omega^\cx,\omega\wedge) \to (B^\cx,\omega\wedge)[1] \to 0$$
defines a morphism
$$a_\omega \colon (B^\cx,\omega\wedge) \to (Z^\cx,\omega\wedge)$$
in $D(X', \mathcal{O}_{X'})$ as well. 
There is also  an inclusion morphism:
$$b_\omega \colon (B^\cx,\omega\wedge) \to (Z^\cx,\omega\wedge).$$
Then it is not difficult to check that
$c_\omega = a_\omega + b_\omega$.  If $a \in k^*$,
we have isomorphisms of complexes 
\begin{eqnarray*}
\lambda_a \colon (Z^\cx,\omega\wedge) & \to & (Z^\cx,a\omega\wedge)\\
\lambda_a \colon (B^\cx,\omega\wedge) & \to & (B^\cx,a\omega\wedge)
\end{eqnarray*}
Using these as identifications, one can check that
$c_{a\omega} = a^{-1} a_\omega + b_\omega$.
This would suggest a negative answer to Question \ref{qomega.q},
but we do not have an example.

\end{remark}

%\bibliography{all}
%\bibliographystyle{plain}

\end{document}